\documentclass[11pt]{amsart}

\usepackage{amsmath,amssymb,amsthm,mathtools}
\usepackage[T1]{fontenc}
\usepackage[utf8]{inputenc}
\usepackage{lmodern}
\usepackage{microtype}
\usepackage{hyperref}

\hypersetup{
  colorlinks=true,
  linkcolor=blue,
  citecolor=blue,
  urlcolor=blue
}

\newtheorem{theorem}{Theorem}[section]
\newtheorem{proposition}[theorem]{Proposition}
\newtheorem{lemma}[theorem]{Lemma}
\newtheorem{corollary}[theorem]{Corollary}

\theoremstyle{remark}
\newtheorem{remark}[theorem]{Remark}

\DeclareMathOperator{\ord}{ord}
\DeclareMathOperator{\CT}{CT}
\DeclareMathOperator{\Span}{Span}
\DeclareMathOperator{\Res}{Res}

\newcommand{\Z}{\mathbf Z}
\newcommand{\Q}{\mathbf Q}
\newcommand{\C}{\mathbf C}
\newcommand{\HH}{\mathfrak H}
\newcommand{\G}{\Gamma_0(3)}

\newcommand{\dd}{\mathrm d}
\newcommand{\Up}{\Lambda_p}
\newcommand{\Vp}{V_p}

\begin{document}

\title[A full $p^{4r}$ supercongruence tower]{A full $p^{4r}$ supercongruence tower for a level-three symmetric-cube hypergeometric sequence}

\author{Alex Shvets}
\address{Haifa, Israel}
\email{alex@shvets.io}
\urladdr{\url{https://shvets.io}}
\thanks{ORCID: 0009-0005-9802-379X}

\begin{abstract}
Let
\[
 {}_2F_1\!\left(\frac13,\frac13;1;27z\right)^3=\sum_{n\ge0}A_nz^n.
\]
We prove the full prime-power supercongruence tower
\[
 A_{mp^r}\equiv A_{mp^{r-1}}\pmod {p^{4r}}
 \qquad (p\ge5,\ m\ge1,\ r\ge1).
\]
The level-three modular expansion used below was already recorded by Moy; the contribution here is the depth-preserving modulus $p^{4r}$, extending the previously known depth-one modulus $p^4$ case to all prime powers.
The proof is modular.  After replacing $A_n$ by $B_n=(-1)^nA_n$, the generating function is realized on $X_0(3)$ by
\[
 \sum_{n\ge0}B_nt^n=\frac{\eta(\tau)^9}{\eta(3\tau)^3},
 \qquad
 t=\frac{\eta(3\tau)^{12}}{\eta(\tau)^{12}}.
\]
The logarithmic derivative
\[
 C=\left(\sum_{n\ge0}B_nt^n\right)\frac{q}{t}\frac{\dd t}{\dd q}
\]
is the Eisenstein series $3E_{5,\chi_0,\chi_3}$.  Lagrange--Buermann gives
\[
 B_m=\CT_q\left(\frac{C(q)}{t(q)^m}\right).
\]
The new point is to replace the one-prime Hecke defect by the prime-power defect
\[
 T_{p^{s+1}}\left(\frac{C}{t^{mp^{s+1}}}\right)
 -T_{p^s}\left(\frac{C}{t^{mp^s}}\right).
\]
Its $i=0$ Hecke layer is the sparse Cartier defect controlling
$A_{mp^{s+1}}-A_{mp^s}$, while all remaining Hecke layers are divisible by the required power of $p$ by induction.  A Fricke involution argument on the two-dimensional space
$M_5(\Gamma_0(3),\chi_3)=\Span\{C,tC\}$ then kills the low-order part exactly and the principal part modulo $p^{4(s+1)}$.
\end{abstract}

\maketitle

\section{Introduction}

Let
\[
 {}_2F_1\!\left(\frac13,\frac13;1;27z\right)^3=\sum_{n\ge0}A_nz^n.
\]
The sequence begins
\[
1,\ 9,\ 135,\ 2439,\ 48519,\ 1023759,\ 22478121,\ldots .
\]
The purpose of this paper is to prove the following prime-power supercongruence tower.

\begin{theorem}\label{thm:main-intro}
For every prime $p\ge5$ and all integers $m,r\ge1$,
\[
 A_{mp^r}\equiv A_{mp^{r-1}}\pmod {p^{4r}}.
\]
\end{theorem}

We separate the new part from the existing level-three dictionary.  Moy \cite[Section 5.2]{Moy} already recorded the exact pair
\[
 f_1=\frac{\eta(\tau)^9}{\eta(3\tau)^3},\qquad
 t=\frac{\eta(3\tau)^{12}}{\eta(\tau)^{12}},
\]
and the expansion
\[
 f_1=1-9t+135t^2-2439t^3+48519t^4-1023759t^5+\cdots .
\]
Moy's framework yields the ordinary $p^r$-level congruences for these level-three expansions; for the first expansion $f_1$ it gives $b_{1,\ell p^r}\equiv b_{1,\ell p^{r-1}}\pmod{p^r}$ for all primes $p\ge5$ (by his Corollary~1.3, which removes the split condition of his Theorem~1.1).  Thus neither the sequence, nor the modular parametrization, nor the ordinary $p^r$ tower is claimed to be new here.  The new point is the depth-preserving strengthening from the ordinary modulus $p^r$ to the modulus $p^{4r}$ for every prime $p\ge5$.

The depth-one case $r=1$, namely
\[
 A_{mp}\equiv A_m\pmod {p^4},
\]
was proved in \cite{ShvetsP4} by a one-prime Hecke--Fricke descent.  The present paper strengthens that result to the full prime-power tower by replacing the one-prime defect with prime-power Hecke lifts and sparse Cartier defects.  A related split-prime supercongruence at the mixed CM point $(1/6,1/3;1)$ was proved in \cite{ShvetsMixed}; that result is not used here.

Put
\[
 q=e^{2\pi i\tau},\qquad
 t(\tau)=\frac{\eta(3\tau)^{12}}{\eta(\tau)^{12}},
 \qquad
 B_n=(-1)^nA_n.
\]
Then
\[
 \sum_{n\ge0}B_nt^n
 ={}_2F_1\!\left(\frac13,\frac13;1;-27t\right)^3
 =\frac{\eta(\tau)^9}{\eta(3\tau)^3}.
\]
Define
\[
 C(q)=\left(\sum_{n\ge0}B_nt(q)^n\right)\frac{q}{t(q)}\frac{\dd t}{\dd q}.
\]
Then $C=3E_{5,\chi_0,\chi_3}$, and hence
\[
 C(q)=1+\sum_{n\ge1}c_nq^n,
 \qquad
 c_n=3\sum_{d\mid n}\chi_3(d)d^4.
\]
The coefficient extraction is given by the Lagrange--Buermann formula
\[
 B_m=\CT_q\left(\frac{C(q)}{t(q)^m}\right).
\]

The one-prime congruence follows from the vanishing modulo $p^4$ of
\[
 \Up\left(\frac{C}{t^{mp}}\right)-\frac{C}{t^m},
\]
where $\Up$ is the Cartier operator $\sum a_nq^n\mapsto \sum a_{np}q^n$.  For the full tower this full $q$-series vanishing is not the right object to iterate.  The correct object is the sparse Cartier defect
\[
 \Up^s\left(
 \Up\left(\frac{C}{t^{mp^{s+1}}}\right)-\frac{C}{t^{mp^s}}
 \right).
\]
The modular lift of this sparse defect is
\[
 \mathcal G_{m,s}:=
 T_{p^{s+1}}\left(\frac{C}{t^{mp^{s+1}}}\right)
 -T_{p^s}\left(\frac{C}{t^{mp^s}}\right).
\]
Prime-power Hecke expansion shows that $\mathcal G_{m,s}$ differs from the sparse Cartier defect by a series divisible by $p^{4(s+1)}$, assuming the lower layers of the tower.  The form $\mathcal G_{m,s}$ is weakly holomorphic of weight $5$ and character $\chi_3$, with poles only at $\infty$.  Since
\[
 M_5(\Gamma_0(3),\chi_3)=\Span\{C,tC\},
\]
it admits a cusp-adapted expansion
\[
 \mathcal G_{m,s}=\beta_{-1}tC+\sum_{j\ge0}\alpha_j\frac{C}{t^j}.
\]
The principal part gives $\alpha_j\equiv0\pmod {p^{4(s+1)}}$ for $j\ge m$.  Applying the Fricke involution $W_3$ gives the exact order bound
\[
 \ord_\infty(\mathcal G_{m,s}|W_3)\ge m+1,
\]
which forces
\[
 \beta_{-1}=0,
 \qquad
 \alpha_j=0\quad(0\le j\le m-1).
\]
Thus $\mathcal G_{m,s}\equiv0\pmod {p^{4(s+1)}}$, and the same is true for the sparse Cartier defect.  Taking constant terms gives Theorem~\ref{thm:main-intro}.

\subsection*{Conventions}
All congruences between Laurent series are coefficientwise congruences in $\Z_{(p)}((q))$.  The operator $\CT_q$ denotes the coefficient of $q^0$ in a Laurent series.  The slash operator used below is determinant-normalized; it is defined explicitly in Section~\ref{sec:fricke}.

\section{The modular parametrization}\label{sec:modular}

Let
\[
 \chi_3(n)=\left(\frac n3\right)
\]
be the quadratic Dirichlet character modulo $3$.
Throughout the paper
\[
 q=e^{2\pi i\tau},\qquad \tau\in\HH.
\]
Set
\[
 t(\tau)=\frac{\eta(3\tau)^{12}}{\eta(\tau)^{12}},
 \qquad
 \Theta(\tau)=\frac{\eta(\tau)^9}{\eta(3\tau)^3}.
\]
Then
\[
 t(q)=q+12q^2+90q^3+508q^4+2391q^5+\cdots
\]
and
\[
 H(q):=\frac{q}{t(q)}=
 \prod_{\substack{n\ge1\\3\nmid n}}(1-q^n)^{12}
 \in 1+q\Z[[q]].
\]

\begin{proposition}[Modular parametrization]\label{prop:modparam}
Let
\[
 {}_2F_1\!\left(\frac13,\frac13;1;27z\right)^3=\sum_{n\ge0}A_nz^n,
 \qquad
 B_n=(-1)^nA_n.
\]
Then
\[
 \sum_{n\ge0}B_nt(\tau)^n=\Theta(\tau).
\]
\end{proposition}

\begin{proof}
Let
\[
 a(q)=\sum_{m,n\in\Z}q^{m^2+mn+n^2},
 \qquad
 b(q)=\frac{\eta(\tau)^3}{\eta(3\tau)},
 \qquad
 c(q)=3\frac{\eta(3\tau)^3}{\eta(\tau)}.
\]
With this normalization, the cubic theory gives
\[
 a(q)^3=b(q)^3+c(q)^3
\]
and
\[
 {}_2F_1\!\left(\frac13,\frac23;1;\frac{c(q)^3}{a(q)^3}\right)=a(q)
\]
\cite[Theorem 2.3 and Corollary 2.4]{BBG}.
Since
\[
 -27t(\tau)=-\frac{c(q)^3}{b(q)^3},
\]
put $z=-27t(\tau)$.  Pfaff's transformation gives
\[
 {}_2F_1\!\left(\frac13,\frac13;1;z\right)
 =(1-z)^{-1/3}
 {}_2F_1\!\left(\frac13,\frac23;1;\frac{z}{z-1}\right).
\]
Here
\[
 \frac{z}{z-1}
 =\frac{c(q)^3}{b(q)^3+c(q)^3}
 =\frac{c(q)^3}{a(q)^3},
\]
and, with the branch normalized to have constant term $1$ at $q=0$,
\[
 (1-z)^{-1/3}
 =\left(1+\frac{c(q)^3}{b(q)^3}\right)^{-1/3}
 =\frac{b(q)}{a(q)}.
\]
Therefore
\[
 {}_2F_1\!\left(\frac13,\frac13;1;-27t(\tau)\right)
 =\frac{b(q)}{a(q)}a(q)
 =\frac{\eta(\tau)^3}{\eta(3\tau)}.
\]
Cubing gives
\[
 \sum_{n\ge0}B_nt(\tau)^n
 ={}_2F_1\!\left(\frac13,\frac13;1;-27t(\tau)\right)^3
 =\frac{\eta(\tau)^9}{\eta(3\tau)^3}.
\]
\end{proof}

Define
\[
 C(\tau)=\Theta(\tau)\frac{Dt(\tau)}{t(\tau)},
 \qquad
 D=q\frac{\dd}{\dd q}=\frac{1}{2\pi i}\frac{\dd}{\dd\tau}.
\]
Equivalently,
\[
 C(q)=\left(\sum_{n\ge0}B_nt(q)^n\right)\frac{q}{t(q)}\frac{\dd t}{\dd q}.
\]

\begin{proposition}[The logarithmic derivative]\label{prop:C-eis}
The form $C$ is the Eisenstein series
\[
 C=3E_{5,\chi_0,\chi_3}.
\]
Consequently
\[
 C(q)=1+\sum_{n\ge1}c_nq^n,
 \qquad
 c_n=3\sum_{d\mid n}\chi_3(d)d^4.
\]
\end{proposition}

\begin{proof}
The eta quotient $\Theta$ has weight $3$ and character $\chi_3$ on $\G$, and $t$ is a Hauptmodul on $X_0(3)$.  Hence $C=\Theta\,Dt/t$ is a meromorphic modular form of weight $5$ and character $\chi_3$ on $\G$.
The $q$-expansions
\[
 \Theta=1-9q+O(q^2),
 \qquad
 t=q+12q^2+O(q^3)
\]
give
\[
 \frac{Dt}{t}=1+12q+O(q^2),
 \qquad
 C=1+3q+O(q^2).
\]

As proved independently in Proposition~\ref{prop:fricke-C}, the form $C$ is also holomorphic at the cusp $0$; hence $C\in M_5(\G,\chi_3)$.
We now recall the needed dimension statement in this special case.  The eta quotient $\Theta$ is nonzero on $\HH$, has order $0$ at $\infty$, and has order $1$ at the cusp $0$ by Proposition~\ref{prop:fricke-C}.  If $f\in S_5(\G,\chi_3)$, then $f/\Theta$ is a holomorphic weight-two form with trivial character on $\G$ and has zero constant term at $\infty$.  Since $X_0(3)$ has genus zero, $S_2(\G)=0$, and $M_2(\G)$ is the one-dimensional Eisenstein space spanned by $E_2(\tau)-3E_2(3\tau)$, whose constant term at $\infty$ is nonzero.  Hence $f/\Theta=0$, and therefore $f=0$.  Thus $S_5(\G,\chi_3)=0$.  The Eisenstein subspace has the two standard Eisenstein series associated with the decompositions $\chi_3=\chi_0\chi_3=\chi_3\chi_0$; equivalently, by the standard dimension formula for modular forms with character on $\Gamma_0(N)$,
\[
 M_5(\G,\chi_3)=\Span\{E_{5,\chi_0,\chi_3},E_{5,\chi_3,\chi_0}\}.
\]
Here
\[
 E_{5,\chi_0,\chi_3}(q)=\frac13+
 \sum_{n\ge1}\left(\sum_{d\mid n}\chi_3(d)d^4\right)q^n
\]
and
\[
 E_{5,\chi_3,\chi_0}(q)=
 \sum_{n\ge1}\left(\sum_{d\mid n}\chi_3(n/d)d^4\right)q^n.
\]
Thus
\[
 3E_{5,\chi_0,\chi_3}=1+3q+O(q^2),
 \qquad
 E_{5,\chi_3,\chi_0}=q+O(q^2).
\]
Since $C$ and $3E_{5,\chi_0,\chi_3}$ have the same constant term and the same coefficient of $q$, their difference is a multiple of $E_{5,\chi_3,\chi_0}$ with zero coefficient of $q$, and therefore is zero.  The divisor-sum formula follows from the displayed Fourier expansion of $E_{5,\chi_0,\chi_3}$.
\end{proof}

\begin{remark}
The same conclusion also follows directly from the usual dimension formula; see, for example, \cite[Sections 3.5--3.6]{DS}.  The quotient argument above is included to make explicit why the cusp-space contribution is zero in this particular weight and level.
\end{remark}

\section{Lagrange--Buermann and integrality}\label{sec:LB}

The next proposition is the coefficient-extraction formula that converts congruences of Laurent series into congruences for $B_m$.

\begin{proposition}[Lagrange--Buermann]\label{prop:LB}
For every $m\ge0$,
\[
 B_m=\CT_q\left(\frac{C(q)}{t(q)^m}\right).
\]
\end{proposition}

\begin{proof}
Since
\[
 \Theta(t)=\sum_{n\ge0}B_nt^n,
\]
one has
\[
 B_m=\CT_t\left(\frac{\Theta(t)}{t^m}\right)
 =\Res_{t=0}\left(\frac{\Theta(t)}{t^m}\frac{\dd t}{t}\right).
\]
The series $t(q)=q+O(q^2)$ is a formal coordinate at $q=0$.  Changing variables gives
\[
 \Res_{t=0}\left(\frac{\Theta(t)}{t^m}\frac{\dd t}{t}\right)
 =\Res_{q=0}\left(\frac{\Theta(t(q))}{t(q)^m}
 \frac{q\,t'(q)}{t(q)}\frac{\dd q}{q}\right).
\]
The numerator is exactly $C(q)$.  Hence
\[
 B_m=\Res_{q=0}\left(\frac{C(q)}{t(q)^m}\frac{\dd q}{q}\right)
 =\CT_q\left(\frac{C(q)}{t(q)^m}\right).
\]
\end{proof}

\begin{corollary}[Integrality]\label{cor:integrality}
For every $m\ge0$, one has $B_m\in\Z$ and hence $A_m\in\Z$.
\end{corollary}

\begin{proof}
By Proposition~\ref{prop:LB},
\[
 B_m=\CT_q\left(\frac{C(q)}{t(q)^m}\right).
\]
Here $C(q)\in\Z[[q]]$ by Proposition~\ref{prop:C-eis}, and
\[
 t(q)^{-m}=q^{-m}H(q)^m,
 \qquad
 H(q)=q/t(q)\in 1+q\Z[[q]].
\]
Thus $C(q)t(q)^{-m}\in q^{-m}\Z[[q]]$, so its constant term is an integer.
Since $B_m=(-1)^mA_m$, this also gives $A_m\in\Z$.
\end{proof}

\begin{remark}[Eisenstein coefficient tower]\label{rem:eis-tower}
The Eisenstein coefficients themselves satisfy the elementary congruence
\[
 c_{mp^r}\equiv c_{mp^{r-1}}\pmod {p^{4r}}
 \qquad(p\ge5,\ m,r\ge1).
\]
Indeed, if $c_n=3\sigma_{4,\chi_3}(n)$ and $m=p^am_0$ with $p\nmid m_0$, then multiplicativity gives
\[
 \sigma_{4,\chi_3}(mp^r)-\sigma_{4,\chi_3}(mp^{r-1})
 =\chi_3(p)^{a+r}p^{4(a+r)}\sigma_{4,\chi_3}(m_0).
\]
This observation is not used in the proof of the full tower; the tower below is obtained from prime-power Hecke defects.
\end{remark}

\begin{remark}
Although the Eisenstein coefficients satisfy the stronger congruence
$c_{mp^r}\equiv c_{mp^{r-1}}\pmod{p^{4r}}$ of Remark~\textup{\ref{rem:eis-tower}},
the formal Moy transfer (Proposition~3.1 of \cite{Moy}) preserves only the
ordinary Dwork level $p^r$ in this setting.  The additional precision $p^{3r}$
is not produced by the formal differential congruence alone; here it is
recovered from the prime-power Hecke expansion (Proposition~\ref{prop:hecke-prime-power})
and the Fricke order drop (Theorem~\ref{thm:sparse-tower}).
\end{remark}

\section{Fricke identities and cusp-adapted expansions}\label{sec:fricke}

We use the determinant-normalized slash operator.  For
\[
 \gamma=\begin{pmatrix}a&b\\c&d\end{pmatrix}\in GL_2^+(\Q)
\]
and a function $f$ of weight $k$, set
\[
 (f|_k\gamma)(\tau)
 =\det(\gamma)^{k/2}(c\tau+d)^{-k}f\left(\frac{a\tau+b}{c\tau+d}\right).
\]
Let
\[
 w_3=\begin{pmatrix}0&-1\\3&0\end{pmatrix},
 \qquad
 W_3f=f|_kw_3
\]
when $f$ has weight $k$.

\begin{proposition}[Fricke identities]\label{prop:fricke-C}
One has
\[
 t|_0W_3=\frac{3^{-6}}{t},
 \qquad
 \Theta|_3W_3=27i\,t\Theta,
 \qquad
 \left(\frac{Dt}{t}\right)\!\Big|_2W_3=-\frac{Dt}{t}.
\]
Consequently
\[
 C|W_3=-27i\,tC,
 \qquad
 (tC)|W_3=-\frac{i}{27}C.
\]
Moreover
\[
 \ord_0(C)=1,
 \qquad
 \ord_0(tC)=0,
 \qquad
 \ord_0\left(\frac{C}{t^j}\right)=j+1\quad(j\ge0).
\]
\end{proposition}

\begin{proof}
We use the standard branch of the square root in the eta transformation formula
\[
 \eta(-1/\tau)=(-i\tau)^{1/2}\eta(\tau).
\]
With this convention the constants below are exactly as written; only their non-vanishing and the resulting orders are used later.
The formula gives immediately
\[
 t\left(-\frac{1}{3\tau}\right)=\frac{3^{-6}}{t(\tau)},
\]
which is the first identity.
For $\Theta$,
\[
 (\Theta|_3W_3)(\tau)=3^{3/2}(3\tau)^{-3}\Theta\left(-\frac{1}{3\tau}\right).
\]
Again using the eta transformation formula,
\[
 \Theta\left(-\frac{1}{3\tau}\right)
 =\frac{\eta(-1/(3\tau))^9}{\eta(-1/\tau)^3}
 =3^{9/2}(-i\tau)^3\frac{\eta(3\tau)^9}{\eta(\tau)^3}.
\]
Thus
\[
 \Theta|_3W_3=27i\frac{\eta(3\tau)^9}{\eta(\tau)^3}=27i\,t\Theta.
\]
For a weight-zero function $f$, the chain rule gives
\[
 D(f|_0\gamma)=(Df)|_2\gamma.
\]
Applying this to $f=t$ and $\gamma=w_3$ gives
\[
 (Dt)|_2W_3=D\left(\frac{3^{-6}}{t}\right)
 =-\frac{3^{-6}}{t^2}Dt.
\]
Dividing by $t|_0W_3=3^{-6}/t$ gives
\[
 \left(\frac{Dt}{t}\right)\!\Big|_2W_3=-\frac{Dt}{t}.
\]
Multiplying with the identity for $\Theta$ yields
\[
 C|W_3=(\Theta|_3W_3)\left((Dt/t)|_2W_3\right)=-27i\,tC.
\]
Since $w_3^2=-3I$, one has $(f|W_3)|W_3=-f$ for weight $5$.  Applying $W_3$ to $C|W_3=-27i\,tC$ gives
\[
 (tC)|W_3=-\frac{i}{27}C.
\]
Finally $\ord_\infty(C)=0$ and $\ord_\infty(tC)=1$, so
\[
 \ord_0(C)=\ord_\infty(C|W_3)=1,
 \qquad
 \ord_0(tC)=\ord_\infty((tC)|W_3)=0.
\]
Since $\ord_0(t)=-1$, the formula
\[
 \ord_0(C/t^j)=j+1
\]
follows.
\end{proof}

\begin{proposition}[Weight-five basis]\label{prop:weight5basis}
One has
\[
 M_5(\G,\chi_3)=\Span\{C,tC\}.
\]
Equivalently,
\[
 C=3E_{5,\chi_0,\chi_3},
 \qquad
 tC=E_{5,\chi_3,\chi_0}.
\]
\end{proposition}

\begin{proof}
The first identity is Proposition~\ref{prop:C-eis}.  The form $tC$ is holomorphic at both cusps by Proposition~\ref{prop:fricke-C}; hence $tC\in M_5(\G,\chi_3)$.  Its expansion begins
\[
 tC=q+15q^2+81q^3+O(q^4),
\]
which is the beginning of $E_{5,\chi_3,\chi_0}$.  Since
\[
 M_5(\G,\chi_3)=\Span\{E_{5,\chi_0,\chi_3},E_{5,\chi_3,\chi_0}\}
\]
and $tC$ has zero constant term and leading coefficient $1$, it follows that
\[
 tC=E_{5,\chi_3,\chi_0}.
\]
The two forms $C$ and $tC$ are linearly independent because their orders at $\infty$ are $0$ and $1$.
\end{proof}

\begin{proposition}[Pole-lowering expansion]\label{prop:polelowering}
Let $f\in M_5^!(\G,\chi_3)$ have poles only at the cusp $\infty$, and suppose the pole order at $\infty$ is at most $N$.  Then
\[
 f=\beta_{-1}tC+\sum_{j=0}^{N}\alpha_j\frac{C}{t^j}
\]
for some $\beta_{-1},\alpha_0,\ldots,\alpha_N\in\C$.
If the $q$-expansion of $f$ is in $\Z_{(p)}((q))$, then all these coefficients lie in $\Z_{(p)}$.
\end{proposition}

\begin{proof}
For $N\ge1$, the form $C/t^N=q^{-N}+O(q^{-N+1})$ has leading coefficient $1$.  Subtracting a unique multiple of $C/t^N$ lowers the pole order by one.  Repeating this process leaves a holomorphic form of weight $5$ and character $\chi_3$, hence an element of $\Span\{C,tC\}$ by Proposition~\ref{prop:weight5basis}.  This gives the expansion.  The case $N=0$ is exactly Proposition~\ref{prop:weight5basis}.

If the $q$-expansion of $f$ has coefficients in $\Z_{(p)}$, then the coefficients in the expansion are recovered by descending triangular elimination.  The diagonal terms are
\[
 \frac{C}{t^j}=q^{-j}+O(q^{-j+1}),
 \qquad
 C=1+O(q),
 \qquad
 tC=q+O(q^2),
\]
whose diagonal entries are $1$.  Since each basis element has $q$-expansion in $\Z[[q,q^{-1}]]$, the recovered coefficients lie in $\Z_{(p)}$.
\end{proof}

\section{Hecke operators and Fricke intertwining}\label{sec:hecke}

For a Laurent series
\[
 f(q)=\sum_{n\gg-\infty}a_nq^n
\]
define
\[
 \Up f=\sum_{n\gg-\infty}a_{np}q^n,
 \qquad
 \Vp f=f(q^p)=\sum_{n\gg-\infty}a_nq^{pn}.
\]
Then $\Up\Vp=1$.
For $p\nmid3$ we use the Hecke normalization compatible with the determinant-normalized slash operator:
\[
 T_pf
 =p^{k/2-1}\left(
 \sum_{b=0}^{p-1}f|_k\begin{pmatrix}1&b\\0&p\end{pmatrix}
 +\chi_3(p)f|_k\begin{pmatrix}p&0\\0&1\end{pmatrix}
 \right)
\]
on weight-$k$ forms with character $\chi_3$.  On $q$-expansions this gives
\[
 T_pf=\Up f+\chi_3(p)p^{k-1}\Vp f.
\]
The prime-power operators $T_{p^a}$ are normalized by the standard Hecke recurrence
\[
 T_pT_{p^a}=T_{p^{a+1}}+\chi_3(p)p^{k-1}T_{p^{a-1}}\qquad(a\ge1),
\]
with $T_1$ the identity; see, for instance, \cite[Chapter 4]{Miyake}.

\begin{proposition}[Prime-power Hecke expansion]\label{prop:hecke-prime-power}
Let $p\ge5$, put $\chi=\chi_3(p)$, and let $a\ge0$.  On $q$-expansions of weakly holomorphic forms of weight $5$ and character $\chi_3$ on $\G$,
\[
 T_{p^a}=
 \sum_{i=0}^{a}\chi^i p^{4i}\Vp^i\Up^{a-i}.
\]
In particular,
\[
 T_p=\Up+\chi p^4\Vp.
\]
\end{proposition}

\begin{proof}
For $a=0$ the formula is the identity operator, and for $a=1$ it is exactly the $q$-expansion formula coming from the normalization fixed above:
\[
 T_p f=\Up f+\chi p^4\Vp f.
\]
Let
\[
 S_a=\sum_{i=0}^{a}\chi^i p^{4i}\Vp^i\Up^{a-i}.
\]
Only the identity $\Up\Vp=1$ is used here; we do not use the false identity $\Vp\Up=1$.  For $i\ge1$ one has $\Up\Vp^i=\Vp^{i-1}$.  Hence
\[
\begin{aligned}
 S_1S_a
 &= (\Up+\chi p^4\Vp)
 \sum_{i=0}^{a}\chi^ip^{4i}\Vp^i\Up^{a-i} \\
 &=\left(\Up^{a+1}+\sum_{j=1}^{a+1}\chi^jp^{4j}\Vp^j\Up^{a+1-j}\right)
 +\chi p^4\sum_{j=0}^{a-1}\chi^jp^{4j}\Vp^j\Up^{a-1-j} \\
 &=S_{a+1}+\chi p^4S_{a-1}\qquad(a\ge1).
\end{aligned}
\]
This is the same recurrence as the Hecke operators
\[
 T_pT_{p^a}=T_{p^{a+1}}+\chi p^4T_{p^{a-1}}
\]
for primes $p\nmid3$ and weight $5$.  Since $S_0=T_1$ and $S_1=T_p$, the displayed formula follows by induction on $a$.
\end{proof}

\begin{lemma}[Holomorphy at the cusp $0$]\label{lem:hecke-cusp0}
If $f\in M_5^!(\G,\chi_3)$ is holomorphic at the cusp $0$, then $T_{p^a}f$ is holomorphic at the cusp $0$ for every $a\ge0$ and every prime $p\ge5$.
\end{lemma}

\begin{proof}
It is enough to prove the assertion for $T_p$, since the prime-power operators are obtained from $T_p$ by the Hecke recurrence.  Write the cusps of $\Gamma_0(3)$ as $\infty$ and $0$.  A reduced fraction $u/v$ is equivalent to $\infty$ precisely when $3\mid v$, and is equivalent to $0$ otherwise.

In the usual formula for $T_p$, the arguments are $p\tau$ and $(\tau+b)/p$ for $0\le b<p$.  If $\tau=u/v$ represents cusp $0$, then $3\nmid v$.  The reduced denominator of $pu/v$ is either $v$ or $v/p$, and the reduced denominator of $(u+bv)/(pv)$ is either $pv$ or $v$.  Since $p\ne3$, all these denominators remain prime to $3$.  Hence the Hecke correspondence over the cusp $0$ maps only to the cusp $0$; equivalently, since $p\nmid3$, the operator $T_p$ preserves the cusp class $0$ on $X_0(3)$.  Therefore no pole at $\infty$ can be pulled back to the cusp $0$, and $T_pf$ is holomorphic at $0$.
\end{proof}

\begin{proposition}[Fricke--Hecke intertwining]\label{prop:fricke-hecke}
For every $a\ge0$ and every prime $p\ge5$,
\[
 T_{p^a}W_3=\chi_3(p)^a W_3T_{p^a}
\]
on $M_5^!(\G,\chi_3)$.
\end{proposition}

\begin{proof}
The case $a=0$ is trivial.  The matrix $w_3$ interchanges the two cusps $0$ and $\infty$; with the determinant-normalized slash operator used here this conjugation does not introduce an additional character factor in the quadratic character case beyond the displayed powers of $\chi_3(p)$.  For $a=1$, use the double-coset representatives
\[
 \alpha_b=\begin{pmatrix}1&b\\0&p\end{pmatrix}
 \quad(0\le b<p),
 \qquad
 \beta=\begin{pmatrix}p&0\\0&1\end{pmatrix}.
\]
With the normalization fixed above, specialized to weight $5$, the Hecke operator is
\[
 T_pf=p^{3/2}\left(\sum_{b=0}^{p-1}f|_5\alpha_b+\chi_3(p)f|_5\beta\right).
\]
For $b\ne0$, choose $b'$ modulo $p$ such that $3bb'\equiv-1\pmod p$, and set
\[
 \gamma_b=
 \begin{pmatrix}
 p&-b'\\
 -3b&(1+3bb')/p
 \end{pmatrix}
 \in\Gamma_0(3).
\]
Then
\[
 w_3\alpha_0=\beta w_3,
 \qquad
 w_3\beta=\alpha_0w_3,
 \qquad
 w_3\alpha_b=\gamma_b\alpha_{b'}w_3\quad(b\ne0).
\]
If $d_b=(1+3bb')/p$ is the lower-right entry of $\gamma_b$, then $d_bp\equiv1\pmod3$, hence
\[
 \chi_3(d_b)=\chi_3(p)^{-1}=\chi_3(p).
\]
Using the transformation law $f|_5\gamma=\chi_3(d)f$ for
$\gamma=\bigl(\begin{smallmatrix}*&*\\*&d\end{smallmatrix}\bigr)\in\Gamma_0(3)$, the preceding identities give
\[
 T_p(f|W_3)=\chi_3(p)(T_pf)|W_3.
\]
This proves the case $a=1$.

For $a\ge2$, the Hecke recurrence
\[
 T_pT_{p^{a-1}}=T_{p^a}+\chi_3(p)p^4T_{p^{a-2}}
\]
and the induction hypothesis give
\[
\begin{aligned}
 T_{p^a}W_3
 &=T_pT_{p^{a-1}}W_3-\chi_3(p)p^4T_{p^{a-2}}W_3 \\
 &=\chi_3(p)^{a-1}T_pW_3T_{p^{a-1}}
   -\chi_3(p)^{a-1}p^4W_3T_{p^{a-2}} \\
 &=\chi_3(p)^aW_3T_pT_{p^{a-1}}
   -\chi_3(p)^{a+1}p^4W_3T_{p^{a-2}} \\
 &=\chi_3(p)^aW_3T_{p^a}.
\end{aligned}
\]
In the third line we used $\chi_3(p)^2=1$.
\end{proof}

\begin{lemma}[Order estimate]\label{lem:order-estimate}
Let $a\ge0$, $p\ge5$, and let
\[
 g(q)=\sum_{n\ge N}a_nq^n
\]
with $N\ge0$.  Then
\[
 \ord_\infty(T_{p^a}g)\ge\left\lceil\frac{N}{p^a}\right\rceil.
\]
\end{lemma}

\begin{proof}
Use Proposition~\ref{prop:hecke-prime-power}.  The term $\Up^{a-i}g$ has order at least $\lceil N/p^{a-i}\rceil$, and applying $\Vp^i$ multiplies the order by $p^i$.  Thus
\[
 \ord_\infty(\Vp^i\Up^{a-i}g)
 \ge p^i\left\lceil\frac{N}{p^{a-i}}\right\rceil
 \ge \left\lceil\frac{N}{p^a}\right\rceil.
\]
Taking the minimum over $0\le i\le a$ gives the claim.
\end{proof}

\section{The sparse Cartier defect}\label{sec:sparse}

For $N\ge0$ put
\[
 f_N=\frac{C}{t^N}.
\]
Since
\[
 C\in\Z[[q]],
 \qquad
 \frac{1}{t}=q^{-1}H(q),
 \qquad
 H(q)\in1+q\Z[[q]],
\]
one has
\[
 f_N=q^{-N}C(q)H(q)^N\in q^{-N}\Z[[q]].
\]

For $s\ge0$ and $m\ge1$, define the sparse Cartier defect
\[
 \mathcal F_{m,s}
 :=\Up^s\left(\Up f_{mp^{s+1}}-f_{mp^s}\right)
 =\Up^{s+1}f_{mp^{s+1}}-\Up^sf_{mp^s}.
\]
Also define the modular Hecke lift
\[
 \mathcal G_{m,s}:=
 T_{p^{s+1}}f_{mp^{s+1}}-T_{p^s}f_{mp^s}.
\]

\begin{lemma}[Pole order of the sparse defect]\label{lem:sparse-pole}
For every $s\ge0$ and $m\ge1$,
\[
 \mathcal F_{m,s}=O(q^{-m+1})
\]
at $\infty$.
\end{lemma}

\begin{proof}
The two Laurent series $\Up f_{mp^{s+1}}$ and $f_{mp^s}$ both begin with $q^{-mp^s}$ with leading coefficient $1$, so their difference is
\[
 \Up f_{mp^{s+1}}-f_{mp^s}=O(q^{-mp^s+1}).
\]
Applying $\Up^s$ keeps only exponents divisible by $p^s$ and divides such exponents by $p^s$.  The smallest multiple of $p^s$ which is at least $-mp^s+1$ is $-(m-1)p^s$.  Hence
\[
 \mathcal F_{m,s}=O(q^{-m+1}).
\]
\end{proof}

\begin{lemma}[Comparison between Hecke and sparse defects]\label{lem:G-F-comparison}
Assume that for all $t<s$ and all $M\ge1$ one has
\[
 \mathcal F_{M,t}\equiv0\pmod {p^{4(t+1)}}.
\]
Then, for every $m\ge1$,
\[
 \mathcal G_{m,s}\equiv \mathcal F_{m,s}\pmod {p^{4(s+1)}}.
\]
For $s=0$ the same conclusion holds without any hypothesis.
\end{lemma}

\begin{proof}
Put $\chi=\chi_3(p)$.  By Proposition~\ref{prop:hecke-prime-power}, the $i=0$ part of $\mathcal G_{m,s}$ is exactly $\mathcal F_{m,s}$.  Therefore
\[
\begin{aligned}
 \mathcal G_{m,s}-\mathcal F_{m,s}
 ={}&\sum_{i=1}^{s}\chi^ip^{4i}\Vp^i
 \left(\Up^{s+1-i}f_{mp^{s+1}}-\Up^{s-i}f_{mp^s}\right)\\
 &+\chi^{s+1}p^{4(s+1)}\Vp^{s+1}f_{mp^{s+1}}.
\end{aligned}
\]
For $1\le i\le s$, the expression in parentheses is
\[
 \Up^{s-i}\left(\Up f_{mp^{s+1}}-f_{mp^s}\right)
 =\mathcal F_{mp^i,s-i},
\]
because $mp^s=(mp^i)p^{s-i}$.  By the induction hypothesis it is divisible by
\[
 p^{4(s-i+1)}.
\]
Multiplication by $p^{4i}$ gives divisibility by $p^{4(s+1)}$.  The final term is visibly divisible by $p^{4(s+1)}$.  This proves the claim.  For $s=0$, the sum over $1\le i\le s$ is empty and only the final term remains.
\end{proof}

\begin{lemma}[Modularity and cusps of $\mathcal G_{m,s}$]\label{lem:G-cusps}
For all $m\ge1$ and $s\ge0$, the form $\mathcal G_{m,s}$ belongs to $M_5^!(\G,\chi_3)$ and has poles only at the cusp $\infty$.
\end{lemma}

\begin{proof}
The form $f_N=C/t^N$ is weakly holomorphic of weight $5$ and character $\chi_3$.  By Proposition~\ref{prop:fricke-C}, it is holomorphic at the cusp $0$.  Hecke operators preserve weak modularity and, by Lemma~\ref{lem:hecke-cusp0}, preserve holomorphy at the cusp $0$.  Thus both terms defining $\mathcal G_{m,s}$ are weakly holomorphic and holomorphic at $0$.  Since $X_0(3)$ has only the two cusps $0$ and $\infty$, all poles of $\mathcal G_{m,s}$ are at $\infty$.
\end{proof}

\section{Proof of the tower}\label{sec:proof}

\begin{theorem}[Sparse Cartier tower]\label{thm:sparse-tower}
For every prime $p\ge5$, every $m\ge1$, and every $s\ge0$,
\[
 \mathcal F_{m,s}\equiv0\pmod {p^{4(s+1)}}.
\]
\end{theorem}

\begin{proof}
We prove the statement by induction on $s$, uniformly in $m$.
Assume first that the statement is known for all lower values $t<s$; when $s=0$ this assumption is empty.
By Lemma~\ref{lem:G-F-comparison},
\[
 \mathcal G_{m,s}\equiv\mathcal F_{m,s}\pmod {p^{4(s+1)}}.
\]
By Lemma~\ref{lem:sparse-pole},
\[
 \mathcal F_{m,s}=O(q^{-m+1}).
\]
Therefore
\[
 [q^{-j}]\mathcal G_{m,s}\equiv0\pmod {p^{4(s+1)}}
 \qquad(j\ge m).
\]

Moreover, $\mathcal G_{m,s}$ has $q$-expansion in $\Z_{(p)}((q))$: this follows from $f_N\in q^{-N}\Z[[q]]$ and the prime-power expansion of Proposition~\ref{prop:hecke-prime-power}.  Therefore, by Lemma~\ref{lem:G-cusps} and Proposition~\ref{prop:polelowering}, there is a finite expansion
\[
 \mathcal G_{m,s}=\beta_{-1}tC+\sum_{j=0}^{J}\alpha_j\frac{C}{t^j}
\]
with coefficients in $\Z_{(p)}$.  The triangular form of the basis at $\infty$ gives, by descending induction on $j$,
\[
 \alpha_j\equiv0\pmod {p^{4(s+1)}}
 \qquad(j\ge m).
\]
Indeed, the coefficient of $q^{-j}$ in $C/t^j$ is $1$, and only the terms $C/t^\ell$ with $\ell\ge j$ can contribute to $q^{-j}$.

It remains to kill the coefficients with $j<m$ and the coefficient of $tC$.  Apply $W_3$.  By Proposition~\ref{prop:fricke-hecke},
\[
 \mathcal G_{m,s}|W_3
 =\chi_3(p)^{s+1}T_{p^{s+1}}(f_{mp^{s+1}}|W_3)
 -\chi_3(p)^sT_{p^s}(f_{mp^s}|W_3).
\]
By Proposition~\ref{prop:fricke-C},
\[
 f_N|W_3=\left(\frac{C}{t^N}\right)\Big|W_3
 =-27i\,3^{6N}Ct^{N+1}.
\]
Thus
\[
 \ord_\infty(f_{mp^{s+1}}|W_3)=mp^{s+1}+1,
 \qquad
 \ord_\infty(f_{mp^s}|W_3)=mp^s+1.
\]
Lemma~\ref{lem:order-estimate} gives
\[
 \ord_\infty T_{p^{s+1}}(f_{mp^{s+1}}|W_3)
 \ge\left\lceil\frac{mp^{s+1}+1}{p^{s+1}}\right\rceil=m+1
\]
and
\[
 \ord_\infty T_{p^s}(f_{mp^s}|W_3)
 \ge\left\lceil\frac{mp^s+1}{p^s}\right\rceil=m+1.
\]
Therefore
\[
 \ord_\infty(\mathcal G_{m,s}|W_3)\ge m+1.
\]

On the other hand,
\[
 (tC)|W_3=-\frac{i}{27}C=-27i\,3^{-6}C
\]
and, for $j\ge0$,
\[
 \left(\frac{C}{t^j}\right)\Big|W_3
 =-27i\,3^{6j}Ct^{j+1}.
\]
Hence
\[
 \mathcal G_{m,s}|W_3
 =-27i\left(\beta_{-1}3^{-6}C+
 \sum_{j=0}^{J}\alpha_j3^{6j}Ct^{j+1}\right).
\]
The forms
\[
 C,\ Ct,\ Ct^2,\ldots
\]
have distinct orders $0,1,2,\ldots$ at $\infty$.  Since $\ord_\infty(\mathcal G_{m,s}|W_3)\ge m+1$, the coefficients of $C,Ct,\ldots,Ct^m$ vanish exactly.  Therefore
\[
 \beta_{-1}=0,
 \qquad
 \alpha_j=0\quad(0\le j\le m-1).
\]
Combining these exact equalities with the congruences for $j\ge m$ gives
\[
 \mathcal G_{m,s}\equiv0\pmod {p^{4(s+1)}}.
\]
Since $\mathcal G_{m,s}\equiv\mathcal F_{m,s}\pmod {p^{4(s+1)}}$, we conclude
\[
 \mathcal F_{m,s}\equiv0\pmod {p^{4(s+1)}}.
\]
This completes the induction on $s$.
\end{proof}

\begin{theorem}[Full $p^{4r}$ tower]\label{thm:main}
Let
\[
 {}_2F_1\!\left(\frac13,\frac13;1;27z\right)^3=\sum_{n\ge0}A_nz^n.
\]
Then, for every prime $p\ge5$ and every pair of integers $m,r\ge1$,
\[
 A_{mp^r}\equiv A_{mp^{r-1}}\pmod {p^{4r}}.
\]
\end{theorem}

\begin{proof}
By Proposition~\ref{prop:LB},
\[
 B_N=\CT_q(f_N).
\]
Let $r=s+1$.  Then
\[
 B_{mp^{s+1}}-B_{mp^s}
 =\CT_q(f_{mp^{s+1}})-\CT_q(f_{mp^s}).
\]
For any Laurent series $g$, $\CT_q(\Up g)=\CT_q(g)$; hence also $\CT_q(\Up^s g)=\CT_q(g)$.  Therefore
\[
 B_{mp^{s+1}}-B_{mp^s}
 =\CT_q\left(\Up^{s+1}f_{mp^{s+1}}-\Up^sf_{mp^s}\right)
 =\CT_q(\mathcal F_{m,s}).
\]
By Theorem~\ref{thm:sparse-tower}, the last constant term is divisible by $p^{4(s+1)}$.  Thus
\[
 B_{mp^r}\equiv B_{mp^{r-1}}\pmod {p^{4r}}.
\]
Finally $B_n=(-1)^nA_n$, and $p$ is odd, so
\[
 (-1)^{mp^r}=(-1)^{mp^{r-1}}.
\]
The same congruence therefore holds for $A_n$.
\end{proof}

\end{document}